\documentclass[nospthms]{svmult}

\usepackage[all]{xy}
\usepackage{amssymb}
\usepackage{amsmath}
\usepackage{amsmath}

\usepackage{epsfig}

\input{psfig.sty}

\newtheorem{proposition}[subsection]{Proposition}

\newtheorem{corollary}[subsection]{Corollary}

\catcode`\Ž=\active \def Ž{\'e}
\catcode`\=\active \def {\`e}
\catcode`\ˆ=\active \def ˆ{\`a} 
\catcode`\=\active \def {\`u}

\catcode`\=\active \def {\^e}        
\catcode`\"=\active \def "{\^{\i}}    
\catcode`\â=\active \def â{\^a}

\catcode`\"=\active \def "{\^{\i}}
\catcode`\ô=\active \def ô{\^o}
\catcode`\ž=\active \def ž{\^u}
\catcode`\é=\active \def é{\`E} 
\catcode`\æ=\active \def æ{\^E}

\def\u{\underline }

\def\oo{\omega}

\def\DD{\Delta}

\def\aa{\alpha}
\def\bb{\beta}

\def\t{\otimes}

\def\End{\mathrm{End}}

\def\KK{{\mathcal K}}

\def\ac{{}^{\scriptstyle \textrm{!`}}}

\def\Ai{A_{\infty}}
\def\AAi{AA_{\infty}}
\def\Ainf#1{AA_{\infty , {#1}}}

\newenvironment{proo}{\begin{trivlist} \item{\emph{Proof.}}}
  {\hfill $\square$ \end{trivlist}}

\def\ZZ{{\mathbb{Z}}}

\def\RR{{\mathbb{R}}}

\def\RR{{\mathbb{R}}}

\def\KK{{\mathbb{K}}}

\def\PP{{\mathcal{P}}}

\def\KKK{{\mathcal{K}}}
\def\TTT{{\mathcal{T}}}
\def\QQQ{{\mathcal{Q}}}

\def\cc{\gamma}

\def\Id{\mathrm{Id }}
\def\id{\mathrm{id }}
\def\Hom{\mathrm{Hom}}

\def\Ksimp#1{\KKK^{#1}_{\rm simp}}
\def\Kcub#1{\KKK^{#1}_{\rm cub}}

\def\epi{\twoheadrightarrow}
\def\mono{\rightarrowtail}

\def\arbreA{\vcenter{\xymatrix@R=3pt@C=3pt{
&& \\
&*{}\ar@{-}[ur] \ar@{-}[ul] \ar@{-}[d]     &\\
&&
}}}

\def\arbreAgrand{\vcenter{\xymatrix@R=30pt@C=30pt{
&& \\
&*{}\ar@{-}[ur] \ar@{-}[ul] \ar@{-}[d]     &\\
&&
}}}

\def\arbreBA{\vcenter{\xymatrix@R=2pt@C=2pt{
&&&&\\
&&&*{}\ar@{-}[ul] & \\
&&*{}\ar@{-}[uurr] \ar@{-}[uull] \ar@{-}[d]     &&\\
&&&&
}}}

\def\arbreAB{\vcenter{\xymatrix@R=2pt@C=2pt{
&&&&\\
&*{}\ar@{-}[ur] &&& \\
&&*{}\ar@{-}[uurr] \ar@{-}[uull] \ar@{-}[d]     &&\\
&&&&
}}}

\def\arbreBB{\vcenter{\xymatrix@R=2pt@C=2pt{
&&*{}&&\\
&&&& \\
&&*{}\ar@{-}[uurr] \ar@{-}[uull] \ar@{-}[d] \ar@{-}[uu]     &&\\
&&&&
}}}

\def\arbreBBbis{\vcenter{\xymatrix@R=4pt@C=4pt@M=0pt{
&&*{}&&\\
&&&& \\
&&*{}\ar@{-}[uurr] \ar@{-}[uull] \ar@{-}[d] \ar@{-}[uu]     &&\\
&&&& }}}

\def\arbreABC{\vcenter{\xymatrix@R=1pt@C=0pt{
&&&&&&\\
&*{}\ar@{-}[ur] &&&&& \\
&&*{}\ar@{-}[uurr] &&&&\\
&&&*{}\ar@{-}[uuurrr] \ar@{-}[uuulll] \ar@{-}[d] &&&\\
&&&&&&
}}}

\def\arbreBAC{\vcenter{\xymatrix@R=1pt@C=0pt{
&&&&&&\\
&&&*{}\ar@{-}[ul] &&& \\
&&*{}\ar@{-}[uurr] &&&&\\
&&&*{}\ar@{-}[uuurrr] \ar@{-}[uuulll] \ar@{-}[d] &&&\\
&&&&&&
}}}

\def\arbreACA{\vcenter{\xymatrix@R=1pt@C=0pt{
&&&&&&\\
&*{}\ar@{-}[ur] &&&&*{}\ar@{-}[ul] & \\
&&&&&&\\
&&&*{}\ar@{-}[uuurrr] \ar@{-}[uuulll] \ar@{-}[d] &&&\\
&&&&&&
}}}

\def\arbreCAB{\vcenter{\xymatrix@R=1pt@C=0pt{
&&&&&&\\
&&&*{}\ar@{-}[ur] &&& \\
&&&&*{}\ar@{-}[uull] &&\\
&&&*{}\ar@{-}[uuurrr] \ar@{-}[uuulll] \ar@{-}[d] &&&\\
&&&&&&
}}}

\def\arbreCBA{\vcenter{\xymatrix@R=1pt@C=0pt{
&&&&&&\\
&&&&&*{}\ar@{-}[ul] & \\
&&&&*{}\ar@{-}[uull] &&\\
&&&*{}\ar@{-}[uuurrr] \ar@{-}[uuulll] \ar@{-}[d] &&&\\
&&&&&&
}}}

\def\arbreAAC{\vcenter{\xymatrix@R=1pt@C=0pt{
&&&&&&\\
&&&&&& \\
&&*{}\ar@{-}[uurr]\ar@{-}[uu]  &&&&\\
&&&*{}\ar@{-}[uuurrr] \ar@{-}[uuulll] \ar@{-}[d] &&&\\
&&&&&&
}}}

\def\arbreCAC{\vcenter{\xymatrix@R=1pt@C=0pt{
&&&&&&\\
&&&*{}\ar@{-}[ul]\ar@{-}[ur] &&& \\
&&&&&&\\
&&&*{}\ar@{-}[uuurrr] \ar@{-}[uuulll]\ar@{-}[uu]  \ar@{-}[d] &&&\\
&&&&&&
}}}

\def\arbreACC{\vcenter{\xymatrix@R=1pt@C=0pt{
&&&&&&\\
&*{}\ar@{-}[ur] &&&&& \\
&&&&&&\\
&&&*{}\ar@{-}[uuurrr] \ar@{-}[uuulll]\ar@{-}[uuu]  \ar@{-}[d] &&&\\
&&&&&&
}}}

\def\arbreCBB{\vcenter{\xymatrix@R=1pt@C=0pt{
&&&&&&\\
&&&&&& \\
&&&&*{}\ar@{-}[uull]\ar@{-}[uu]  &&\\
&&&*{}\ar@{-}[uuurrr] \ar@{-}[uuulll] \ar@{-}[d] &&&\\
&&&&&&
}}}

\def\arbreCCA{\vcenter{\xymatrix@R=1pt@C=0pt{
&&&&&&\\
&&&&&*{}\ar@{-}[ul] & \\
&&&&&&\\
&&&*{}\ar@{-}[uuurrr] \ar@{-}[uuulll]\ar@{-}[uuu]  \ar@{-}[d] &&&\\
&&&&&&
}}}

\def\arbreBBC{\vcenter{\xymatrix@R=1pt@C=0pt{
&&&&&&\\
&&&&&& \\
&&*{}\ar@{-}[uu] \ar@{-}[uurr] &&&&\\
&&&*{}\ar@{-}[uuurrr] \ar@{-}[uuulll] \ar@{-}[d] &&&\\
&&&&&&
}}}

\def\arbreCCC{\vcenter{\xymatrix@R=1pt@C=0pt{
&&&&&&\\
&&&&&& \\
&&&&&&\\
&&&*{}\ar@{-}[uuurrr]  \ar@{-}[uuulll]\ar@{-}[uuul] \ar@{-}[uuur]  \ar@{-}[d] &&&\\
&&&&&&
}}}

\def\arbreCCCbis{\vcenter{\xymatrix@R=4pt@C=4pt@M=0pt{
&&&&&&\\
&&&&&& \\
&&&&&&\\
&&&*{}\ar@{-}[uuurrr]  \ar@{-}[uuulll]\ar@{-}[uuul] \ar@{-}[uuur]  \ar@{-}[d] &&&\\
&&&&&& }}}

\def\arbreABCD{\vcenter{\xymatrix@R=1pt@C=0pt{
&&&&&&&&\\
&*{}\ar@{-}[ur] &&&&&&& \\
&&*{}\ar@{-}[uurr] &&&&&&\\
&&&*{}\ar@{-}[uuurrr] &&&&&\\
&&&&*{}\ar@{-}[uuuurrrr] \ar@{-}[uuuullll] \ar@{-}[d] &&&&\\
&&&&&&&&
}}}

\def\arbreBACD{\vcenter{\xymatrix@R=1pt@C=0pt{
&&&&&&&&\\
&&&*{}\ar@{-}[ul] &&&&& \\
&&*{}\ar@{-}[uurr] &&&&&&\\
&&&*{}\ar@{-}[uuurrr] &&&&&\\
&&&&*{}\ar@{-}[uuuurrrr] \ar@{-}[uuuullll] \ar@{-}[d] &&&&\\
&&&&&&&&
}}}

\def\arbreCABD{\vcenter{\xymatrix@R=1pt@C=0pt{
&&&&&&&&\\
&&&*{}\ar@{-}[ur] &&&&& \\
&&&&*{}\ar@{-}[uull] &&&&\\
&&&*{}\ar@{-}[uuurrr] &&&&&\\
&&&&*{}\ar@{-}[uuuurrrr] \ar@{-}[uuuullll] \ar@{-}[d] &&&&\\
&&&&&&&&
}}}

\def\arbreDABC{\vcenter{\xymatrix@R=1pt@C=0pt{
&&&&&&&&\\
&&&*{}\ar@{-}[ur]& &&&& \\
&&&&*{}\ar@{-}[uurr] &&&&\\
&&&&&*{}\ar@{-}[uuulll] &&&\\
&&&&*{}\ar@{-}[uuuurrrr] \ar@{-}[uuuullll] \ar@{-}[d] &&&&\\
&&&&&&&&
}}}

\def\arbreDBAC{\vcenter{\xymatrix@R=1pt@C=0pt{
&&&&&&&&\\
&&&& &*{}\ar@{-}[ul]&&& \\
&&&&*{}\ar@{-}[uurr] &&&&\\
&&&&&*{}\ar@{-}[uuulll] &&&\\
&&&&*{}\ar@{-}[uuuurrrr] \ar@{-}[uuuullll] \ar@{-}[d] &&&&\\
&&&&&&&&
}}}

\def\arbreDCAB{\vcenter{\xymatrix@R=1pt@C=0pt{
&&&&&&&&\\
&&&& &*{}\ar@{-}[ur]&&& \\
&&&& &&*{}\ar@{-}[uull] &&\\
&&&&&*{}\ar@{-}[uuulll] &&&\\
&&&&*{}\ar@{-}[uuuurrrr] \ar@{-}[uuuullll] \ar@{-}[d] &&&&\\
&&&&&&&&
}}}

\def\arbreDCBA{\vcenter{\xymatrix@R=1pt@C=0pt{
&&&&&&&&\\
&&&& &&&*{}\ar@{-}[ul]& \\
&&&& &&*{}\ar@{-}[uull] &&\\
&&&&&*{}\ar@{-}[uuulll] &&&\\
&&&&*{}\ar@{-}[uuuurrrr] \ar@{-}[uuuullll] \ar@{-}[d] &&&&\\
&&&&&&&&
}}}

\def\arbreBBDD{\vcenter{\xymatrix@R=1pt@C=0pt{
&&&&&&&&\\
& &&&&&&& \\
&&*{}\ar@{-}[uurr] \ar@{-}[uu]&&&&&&\\
&&&&&&&&\\
&&&&*{}\ar@{-}[uuuurrrr]\ar@{-}[uuuur] \ar@{-}[uuuullll] \ar@{-}[d] &&&&\\
&&&&&&&&
}}}

\def\arbreBADA{\vcenter{\xymatrix@R=1pt@C=0pt{
&&&&&&&&\\
&&&*{} \ar@{-}[ul]&&&&*{}\ar@{-}[ul] & \\
&&*{}\ar@{-}[uurr]&&&&&&\\
&&&&&&&&\\
&&&&*{}\ar@{-}[uuuurrrr]\ar@{-}[uuuullll] \ar@{-}[d] &&&&\\
&&&&&&&&
}}}

\def\arbreADBA{\vcenter{\xymatrix@R=1pt@C=0pt{
&&&&&&&&\\
&*{}\ar@{-}[ur]&&&&&&*{}\ar@{-}[ul] & \\
&&&&&&*{}\ar@{-}[uull]&&\\
&&&&&&&&\\
&&&&*{}\ar@{-}[uuuurrrr]\ar@{-}[uuuullll] \ar@{-}[d] &&&&\\
&&&&&&&&
}}}

\def\arbreCCCD{\vcenter{\xymatrix@R=1pt@C=0pt{
&&&&&&&&\\
& &&&&&& & \\
&&&&&& &&\\
&&&*{}\ar@{-}[uuur]\ar@{-}[uuul]\ar@{-}[uuurrr]&&&&&\\
&&&&*{}\ar@{-}[uuuurrrr]\ar@{-}[uuuullll] \ar@{-}[d] &&&&\\
&&&&&&&&
}}}

\def\arbreDADD{\vcenter{\xymatrix@R=1pt@C=0pt{
&&&&&&&&\\
& &&&*{}\ar@{-}[ur]\ar@{-}[ul]&&& & \\
&&&&&& &&\\
&&&&&&&&\\
&&&&*{}\ar@{-}[uuuurrrr]\ar@{-}[uuuullll] \ar@{-}[d]\ar@{-}[uuuurr]\ar@{-}[uuu] &&&&\\
&&&&&&&&
}}}

\def\arbreDBBD{\vcenter{\xymatrix@R=1pt@C=0pt{
&&&&&&&&\\
& &&&&&& & \\
&&&&*{}\ar@{-}[uurr]\ar@{-}[uull]\ar@{-}[uu]&& &&\\
&&&&&&&&\\
&&&&*{}\ar@{-}[uuuurrrr]\ar@{-}[uuuullll] \ar@{-}[d]\ar@{-}[uu] &&&&\\
&&&&&&&&
}}}

\def\arbreDDAD{\vcenter{\xymatrix@R=1pt@C=0pt{
&&&&&&&&\\
& &&&&*{}\ar@{-}[ur]\ar@{-}[ul]&& & \\
&&&&&& &&\\
&&&&&&&&\\
&&&&*{}\ar@{-}[uuuurrrr]\ar@{-}[uuuullll] \ar@{-}[d]\ar@{-}[uuuull]\ar@{-}[uuur] &&&&\\
&&&&&&&&
}}}

\def\arbreDDBB{\vcenter{\xymatrix@R=1pt@C=0pt{
&&&&&&&&\\
& &&&&& & & \\
&&&&&&*{}\ar@{-}[uu]\ar@{-}[uull] &&\\
&&&&&&&&\\
&&&&*{}\ar@{-}[uuuurrrr]\ar@{-}[uuuullll] \ar@{-}[d]\ar@{-}[uuuull] &&&&\\
&&&&&&&&
}}}

\def\KzeroF{\xymatrix@R=4pt@C=4pt{
\\
\\
\\
\bullet\\
}}

\def\KunF{\xymatrix@R=4pt@C=4pt{
&&\\
&&\\
&&\\
*{}\ar@{->}[rrr]&&&*{}\\
}}

\def\KdeuxF{\xymatrix@R=4pt@C=4pt{
&& &&\\
&&*{}\ar[dll]\ar[ddrrr]  &&&&\\
 *{}\ar[dd] &&&&&\\
&&&&&*{}\ar[ddlll] \\
 *{}\ar[drr] &&&&& \\
&&*{}&&& \\
&& &&&
}}

\def\KdeuxCubical{\xymatrix@R=4pt@C=4pt{
&&& &&\\
&&&*{}\ar@{=>}[dll]\ar@{=>}[ddrrr]  &&&&\\
& *{}\ar@{=}[dd] &&&&&\\
\qquad \mapsto \qquad &&&&&&*{}\ar@{=>}[ddlll] \\
& *{}\ar@{=>}[drr] &&&&& \\
&&&*{}&&& \\
&&& &&&
}}

\def\KdeuxSimp{\vcenter{\xymatrix@R=4pt@C=4pt{
&& &&\\
&&*{}\ar[dll]\ar[ddrrr]\ar@{=>}[dddd]   &&&&\\
 *{}\ar[dd]\ar[dddrr]  &&&&&\\
&&&&&*{}\ar[ddlll] \\
 *{}\ar[drr] &&&&& \\
&&*{}&&& \\
&& &&&
}}}

\def\KdeuxPushed{\vcenter{\xymatrix@R=4pt@C=4pt{
&& &&\\
&&*{}\ar@{=}[dll]\ar[ddrrr]  &&&&\\
 *{}\ar@{=}[dd] &&&&&\\
&&&&&*{}\ar[ddlll] \\
 *{}\ar@{=>}[drr] &&&&& \\
&&*{}&&& \\
&& &&&
}}}

\def\KtroisF{\xymatrix@R=0pt@C=0pt{
&*{}\ar[rrrrrrrrrr] *{}\ar@{<-}[rrrrddd] *{}\ar[ldd] &&&& &&& &&&*{}\ar@{.>}[llldddddd]  *{}\ar@{<-}[rrd] && \\
&&&& &&& &&& &&&*{}\ar@{<-}[llldd] *{}\ar[llldddddd]   \\
*{}\ar@{<-}[rrrrddd]  *{}\ar[rrdddd] &&&&& &&& &&& && \\
&&&&&*{}\ar[rrrrr] *{}\ar[ldd]  &&& &&*{}\ar[lldddd] & && \\
&&&&& &&& &&& && \\
&&&&*{}\ar[rdd] & &&& &&& & &\\
&&*{}\ar@{.>}[rrrrrr] *{}\ar@{<-}[rrrddd] &&& &&&*{}\ar@{<.}[rrd]  &&& && \\
&&&&&*{}\ar[rrr] *{}\ar[dd]  &&&*{}\ar[dd]  &&*{}\ar@{<-}[lldd] & & &\\
&&&&& &&& &&& && \\
&&&&&*{}\ar[rrr]  &&&*{} &&& && \\
}}

\def\KdeuxTfat{\xymatrix@R=10pt@C=10pt{
&&*{}\ar@{=>}[dll]\ar@{=>}[ddrrr]\ar[dddd] &&&&\\
 *{}\ar@{=>}[dd] \ar[dddrr]&b &&&&\\
&&& a&&*{}\ar[ddlll] \\
 *{}\ar[drr] &c&&&& \\
&&*{}&&& \\
}}

\def\KdeuxT{\xymatrix@R=10pt@C=10pt{
&&&&&&\\
&&&&&&\\
&&*{}\ar[dll]\ar[ddrrr]\ar[dddd] &&&&\\
 *{}\ar[dd] \ar[dddrr]& &&&&\\
&&& &&*{}\ar[ddlll] \\
 *{}\ar[drr] &&&&& \\
&&*{}&&& \\
}}

\def\KdeuxMainS{\xymatrix@R=10pt@C=10pt{
&&*{}\ar[dll]\ar@{=>}[ddrrr]\ar@{=>}[dddd] &&&&\\
 *{}\ar[dd] \ar[dddrr]& &&&&\\
&&& &&*{}\ar@{=>}[ddlll] \\
 *{}\ar[drr] &&&&& \\
&&*{}&&& \\
}}

\def\KunSquare{\xymatrix@R=20pt@C=20pt{
&&&&\\
*{}\ar@{|-}[r]&*{}\ar@{|-|}[r]&&&\\
}}

\def\KdeuxSquare{\xymatrix@R=6pt@C=6pt{
& &*{}\ar@{-}[dl] \ar@{-}[dr]& &  &\\
&*{}\ar@{-}[dl]\ar@{-}[ddr] & &*{}\ar@{-}[ddl]\ar@{-}[ddrr]] &  &\\
*{}\ar@{-}[d]& & & &  &\\
*{}\ar@{-}[d] \ar@{-}[rr]& &*{}\ar@{-}[ddl] \ar@{-}[ddr]& &  &*{}\ar@{-}[ddll]\\
*{}\ar@{-}[dr]& & & &  &\\
& *{}\ar@{-}[dr]& & *{}\ar@{-}[dl] &  &\\
& & *{}& &  &\\
}}

\def\KdeuxSquareMain{\xymatrix@R=6pt@C=6pt{
& &*{}\ar@{=>}[dl] \ar@{=>}[dr]& &  &\\
&*{}\ar@{-}[dl]\ar@{=>}[ddr] & &*{}\ar@{=>}[ddl]\ar@{-}[ddrr]] &  &\\
*{}\ar@{-}[d]& & & &  &\\
*{}\ar@{-}[d] \ar@{-}[rr]& &*{}\ar@{-}[ddl] \ar@{-}[ddr]& &  &*{}\ar@{-}[ddll]\\
*{}\ar@{-}[dr]& & & &  &\\
& *{}\ar@{-}[dr]& & *{}\ar@{-}[dl] &  &\\
& & *{}& &  &\\
}}
\def\KzeroTree{\xymatrix@R=4pt@C=4pt{
&& \\
&*{}\ar@{-}[ur] \ar@{-}[ul] \ar@{-}[d]     &\\
&& \\
&\bullet&\\
&1&
}}

\def\KunTree{\xymatrix@R=2pt@C=2pt{
&&&&    &&   &&&&\\
&*{}\ar@{-}[ur] &&&     &&    &&&*{}\ar@{-}[ul] & \\
&&*{}\ar@{-}[uurr] \ar@{-}[uull] \ar@{-}[d]     &&    &&    &&*{}\ar@{-}[uurr] \ar@{-}[uull] \ar@{-}[d]     &&\\
&&&&    &&    &&&& \\
&&*{}\ar[rrrrrr] &&    &&    &&&&\\
&&12 &&&&&&21&&
}}

\def\KdeuxTree{\xymatrix@R=4pt@C=10pt{
& & 123 & & &  \\
& & *{}\ar[dl] \ar[ddrr] & & & \\
213 & *{}\ar[dd] & & & & \\
 & & & & *{}\ar[ddll] & 141 \\
312 &*{}\ar[dr] & & & & \\
 & & *{} & & & \\
  & & 321 & & & \\
}}

\def\KtroisCube{\xymatrix@R=10pt@C=10pt{
& & *{}\ar[rrrr] && & &*{}  \\
 & & && &*{}\ar[ur] & \\
*{}\ar[rrrr]\ar[uurr]  &&\ar[uu]  &&*{}\ar[ur] & & \\
 & & *{}\ar@{-}[u]\ar@{-}[rr] &&\ar@{-}[r]  &\ar[r]  &  *{}\ar[uuu] \\
  & & && & & *{}\ar[u] \\
   & & && &*{}\ar[ur]\ar[uuuu] & \\
    & &*{}\ar[uuu]\ar@{-}[rr] &&\ar@{-}[r]  &\ar[r]  &*{}\ar[uu]  \\
     & & && &*{}\ar[ur]\ar[uu]  & \\
*{}\ar[rrrr]\ar[uurr] \ar[uuuuuu]  & & &&*{}\ar[ur]\ar[uuuuuu]  & & \\
 }}

\def\pforTriangleUn{\xymatrix@R=10pt@C=10pt{
&&\bullet\ar[dll]\ar[ddrrr]\ar[dddd] &&&&  &\\
 *{}\ar[dd] \ar[dddrr]& &&&&  &\\
&&& &&\bullet\ar[ddlll]  &  \epi  \\
 *{}\ar[drr] &&&&&  & \\
&&\bullet&&&  & \\
}}

\def\pforTriangleDeux{\xymatrix@R=10pt@C=10pt{
\bullet\ar[ddrrr]\ar@{=>}[dddd] &&&&  &\\
&&&  &\\
& &&\bullet\ar[ddlll]  & \mono \\
&&&  &\\
\bullet&&&  &\\
}}

\def\pforTriangleTrois{\xymatrix@R=10pt@C=10pt{
&&\bullet\ar@{-}[dll]\ar[ddrrr] &&&& \\
 *{}\ar@{-}[dd] & &&&&  &\\
&&& &&\bullet\ar[ddlll]  \\
 *{}\ar[drr] &&&&&  \\
&&\bullet&&& \\
}}

\def\pforCubeUn{\xymatrix@R=10pt@C=10pt{
& &\bullet\ar[dl] \ar[dr]& &  &  &\\
&\bullet\ar[ddr] & &\bullet\ar[ddl]&  &  &\\
*{}\ar[ur]\ar[d]& & & &  &  &\\
*{}\ar[rr]& &\bullet& &  &*{}\ar[ddll]\ar[uull]   & \epi  \\
*{}\ar[dr]\ar[u] & & & &  &  &\\
& *{}\ar[uur] & & *{}\ar[uul]&  &  &\\
& & *{}\ar[ul]\ar[ur] & &  &  &\\
}}
\def\pforCubeDeux{\xymatrix@R=10pt@C=10pt{
&& &&&  &\\
&&\bullet\ar[ddll]\ar[ddrr]&&&  &\\
&&&&  &\\
\bullet\ar@{=>}[ddrr] && &&\bullet\ar@{=>}[ddll]  & \mono \\
&&&&  &\\
&&\bullet&&  &\\
}}

\def\pforCubeTrois{\xymatrix@R=10pt@C=10pt{
&& &&&  &\\
&&\bullet\ar[dll]\ar[ddrrr] &&&& \\
 \bullet\ar@{-}[dd] & &&&&  &\\
&&& &&\bullet\ar[ddlll]  \\
 *{}\ar[drr] &&&&&  \\
&&\bullet&&& \\
}}

\def\KcubesimpUn{\xymatrix@R=10pt@C=10pt{
&&&&&&\\
& &*{}\ar[dl] \ar[dr]& &  &  &\\
&*{}\ar[ddr] & &*{}\ar[ddl]&  &  &\\
*{}\ar[ur]\ar[d]& & & &  &  &\\
*{}\ar[rr]& &*{}& &  &*{}\ar[ddll]\ar[uull]   & \to  \\
*{}\ar[dr]\ar[u] & & & &  &  &\\
& *{}\ar[uur] & & *{}\ar[uul]&  &  &\\
& & *{}\ar[ul]\ar[ur] & &  &  &\\
}}

\def\KcubesimpDeuxbis{\xymatrix@R=10pt@C=10pt{
                        &                        &    &*{}\ar[dll] \ar[ddrr]\ar[dddd]&   &                             &&  &\\
                        &*{}\ar[dddrr]& &                              &  &&                     &  &\\
*{}\ar[ur]\ar[dd]\ar[ddrrr]&           &  &                            &   & *{}\ar[ddll]&                                         &  &\\
                        &                    &     &                               &    & &                                          &  &\\
*{}\ar[rrr]           &                      &     &*{}                   &    &  &                      &*{}\ar[ddll]\ar[uull]\ar[llll]   & \leftarrow \\
                           &                     &     &                       &  &   &                     &                                            &   \\
*{}\ar[dr]\ar[uu]\ar[uurrr] &            &     &                              &      &                 *{}\ar[uull]                      &  &  &\\
                        & *{}\ar[uuurr]     &     &                             &    &  &                          &  &\\
                        &                      &      &*{}\ar[ull]\ar[uurr]\ar[uuuu]&     &                           &&  &\\
}}

\def\KcubesimpTrois{\xymatrix@R=10pt@C=10pt{
&&&&&&&\\
&&&&&&&\\
&&*{}\ar[dll]\ar[ddrrr]\ar[dddd] &&&&  &\\
 *{}\ar[dd] \ar@{=>}[dddrr]& b&&&&  &\\
&&&c &&*{}\ar@{=>}[ddlll]  &  \\
 *{}\ar@{=>}[drr] &a&&&&  & \\
&&*{}&&&  & 
}}

\def\KtroisTfat{\xymatrix@R=6pt@C=6pt{
&*{}\ar@{=>}[rrrrrrrrrr] *{}\ar@{<=}[rrrrddd] *{}\ar@{=>}[ldd] &&&& &&& &&&*{}\ar@{<=}[rrd] && \\
&&&& &&& &&& &&&*{}\ar@{<=}[llldd] *{}\ar@{=>}[lldddddd]   \\
*{}\ar@{<=}[rrrrddd]  *{}\ar@{=>}[rrdddd] &&&&& &&& &&& && \\
&&&&&*{}\ar@{=>}[rrrrr] \ar@{=>}[ldd] \ar[ulllll] \ar[uuurrrrrr] \ar[ddddrrr] &&& &&*{}\ar@{=>}[lldddd] \ar[uuur] \ar[rdddd]& && \\
&&&&& && & &&& && \\
&&&&*{}\ar@{=>}[rdd] \ar[dll] \ar[ddrrrr]& &&& &&& & &\\
&&*{}\ar@{<=}[rrrddd] &&& &&&*{}&&& && \\
&&&&&*{}\ar@{=>}[rrr] *{}\ar@{=>}[dd] \ar[ulll] \ar[ddrrr]&&&*{}\ar@{=>}[dd] \ar[rrr] &&&*{}\ar@{<=}[llldd]  & &\\
&&&&& &&& &&& && \\
&&&&&*{}\ar@{=>}[rrr]  &&&*{} &&& && \\
}}

\def\KtroisMainS{\xymatrix@R=6pt@C=6pt{
&*{}\ar@{=>}[rrrrrrrrrr] *{}\ar@{<=}[rrrrddd] *{}\ar[ldd] &&&& &&& &&&*{}&& \\
&&&& &&& &&& &&&*{}\ar@{<-}[llldd] *{}\ar[lldddddd]\ar[llu]     \\
*{}\ar@{<-}[rrrrddd]\ar[rrdddd] &&&&& &&& &&& && \\
&&&&&*{}\ar[rrrrr] \ar[ldd] \ar[ulllll] \ar@{=>}[uuurrrrrr] \ar[ddddrrr] &&& &&*{}\ar[lldddd] \ar[uuur] \ar[rdddd]& && \\
&&&&& && & *{}\ar@{=>}[ulll]\ar@{=>}[uuuulllllll]\ar@{.}[ddllllll] \ar@{.}[dddrrr]
\ar@{<=}[uuuurrr]  &&& && \\
&&&&*{}\ar[rdd] \ar[dll] \ar[ddrrrr]& &&&&&& & &\\
&&*{}\ar[rrrddd] &&& &&&*{}&&& && \\
&&&&&*{}\ar[rrr] *{}\ar[dd] \ar[ulll] \ar[ddrrr]&&&*{}\ar[dd] \ar[rrr] &&&*{} & &\\
&&&&& &&& &&& && \\
&&&&&*{}\ar[rrr]  &&&*{} \ar[rrruu] &&& && \\
}}

\def\KtroisDeformed{\xymatrix@R=6pt@C=6pt{
&*{}\ar@{=>}[rrrrrrrrrr] *{}\ar@{<=}[rrrrddd] *{}\ar@{=}[ldd] &&&& &&& &&&*{}&& \\
&&&& &&& &&& &&&*{}\ar@{=}[llldd] *{}\ar[lldddddd] \ar@{=>}[llu]   \\
*{}\ar@{<-}[rrrrddd]  *{}\ar@{=}[rrdddd] &&&&& &&& &&& && \\
&&&&&*{}\ar@{=}[rrrrr] \ar@{=}[ldd]  \ar@[uuurrrrrr]  &&& &&*{}\ar[lldddd] & && \\
&&&&& &&
*{} \ar@{<=}[ddlllll] \ar@{<=}[dddrrrr] \ar@{<=}[uuuurrrr] & &&& && \\
&&&&*{}\ar@{=}[rdd] & &&& &&& & &\\
&&*{}\ar[rrrddd] &&& &&&*{}&&& && \\
&&&&&*{}\ar@{=}[rrr] *{}\ar[dd] &&&*{}\ar@{=}[dd] &&&*{}& &\\
&&&&& &&& &&& && \\
&&&&&*{}\ar[rrr]  &&&*{}\ar@{=}[rrruu]   &&& && \\
}}

\def\KdeuxDIAG{\xymatrix@R=10pt@C=10pt{
& &*{}\ar@{-}[ddll]\ar@{-}[ddddrrrr] &    & & & \\
&*{}\ar@{-}[dd] & &    & & & \\
*{}\ar@{-}[dddd]& & &    &*{}\ar@{-}[ddll] & & \\
&*{}\ar@{-}[dl]\ar@{-}[dr] & &    & & & \\
*{}\ar@{-}[dr]& &*{}\ar@{-}[dl]\ar@{-}[ddrr] &    & & &*{}\ar@{-}[ddddllll] \\
& *{}\ar@{-}[dd]& &    & & & \\
*{}\ar@{-}[ddrr]& & &    &*{} & & \\
&*{} &*{} &    & & & \\
& &  *{} &  & & & 
}}

\begin{document}

\title*{The diagonal of the Stasheff polytope}
\author{Jean-Louis Loday
}
\institute{Institut de Recherche Math\'ematique Avanc\'ee\\
    CNRS et Universit\'e Louis Pasteur\\
    7 rue R.~Descartes\\
    67084 Strasbourg Cedex, France\\
\texttt{loday@math.u-strasbg.fr}}
%
%
\maketitle
\centerline{\emph{To Murray Gerstenhaber and Jim Stasheff}}

\vskip1cm

\date{\today}

\begin{abstract} We construct an A-infinity structure on the tensor product of two A-infinity algebras by using the simplicial decomposition of the Stasheff polytope. The key point is the construction of an operad AA-infinity based on the simplicial Stasheff polytope. The operad AA-infinity admits a coassociative diagonal  and the operad A-infinity is a retract by deformation of it.  We compare these constructions with analogous constructions due to Saneblidze-Umble and Markl-Shnider based on the Boardman-Vogt cubical decomposition  of the Stasheff polytope.
\keywords{Stasheff polytope, associahedron, operad, bar-cobar construction, cobar construction, A-infinity algebra, AA-infinity algebra, diagonal.}
\end{abstract}

\section*{Introduction} \label{S:int}

An associative algebra up to homotopy, or $A_{\infty}$-algebra,  is a chain complex $(A,d_{A})$ equipped with an $n$-ary operation $\mu_{n}$ for each $n\geq 2$ verifying $\mu \circ \mu =0$. See \cite{Stasheff}, or, for instance,  \cite{Keller}. Here we put
$$\mu := d_{A} + \mu_{2} + \mu_{3}+\cdots\    : T(A) \to T(A),$$
 where $\mu_{n}$ has been extended to the tensor coalgebra $T(A)$ by coderivation. In particular $\mu_{2}$ is not associative, but only associative up to homotopy in the following sense:
$$\mu_{2}\circ (\mu_{2}\t \id) - \mu_{2}\circ (\id\t \mu_{2}) = d_{A}\circ \mu_{3}+ \mu_{3}\circ d_{A^{\t 3}}\ .$$

Putting an $A_{\infty}$-algebra structure on the tensor product of two $A_{\infty}$-algebras is a long standing problem, cf.~for instance \cite{Proute, GZ}. Recently a solution has been constructed by Saneblidze and Umble, cf.~\cite{SU1,SU2}, by constructing a diagonal $\Ai \to \Ai\t \Ai$ on the operad $\Ai$ which governs the $\Ai$-algebras. Recall that, over  a field,  the operad $\Ai$ is the minimal model of the operad $As$ governing the associative algebras. The differential graded module $(\Ai)_{n}$ of the $n$-ary operations is the chain complex of the Stasheff polytope. In \cite{MS} Markl and Shnider give a conceptual construction of the Saneblidze-Umble diagonal by using the Boardman-Vogt model of $As$. This model is the bar-cobar construction on $As$, denoted $\Omega B As$, in the operadic framework. It turns out that there exists a coassociative diagonal on $\Omega B As$. This diagonal, together with the quasi-isomorphisms $q: \Ai {\to} \Omega B As$ and  $p:\Omega B As{\to}  \Ai$ permit them to construct a diagonal on $\Ai$ by composition:
$$ \Ai \stackrel{q}{\to} \Omega B As \to \Omega B As \t \Omega B As \stackrel{p\t p}{\longrightarrow} \Ai \t \Ai\ .$$

The aim of this paper is to give an alternative solution to the diagonal problem by relying on the \emph{simplicial decomposition of the Stasheff polytope} described in \cite{JLLparking}. It leads to a new   model $\AAi  $ of the operad $As$, whose dg module $(\AAi)_{n}$ is the chain complex of a simplicial decomposition of the Stasheff polytope.   Because of its simplicial nature, the operad $\AAi$ has a coassociative diagonal (by means of the Alexander-Whitney map) and therefore we get a new diagonal on $\Ai$ by composition:
$$ \Ai \stackrel{q'}{\to} \AAi \to \AAi \t \AAi \stackrel{p'\t p'}{\longrightarrow} \Ai \t \Ai\ .$$

The map $q': \Ai \to \AAi$ is induced by the simplicial decomposition of the associahedron. The map $p':\AAi\to \Ai$ is slightly more involved to construct. It is induced by the deformation of the ``main simplex'' of the associahedron into the big cell of the associahedron. Here the main simplex is defined by the shortest path in the Tamari poset structure of the planar binary trees.

We compute the diagonal map on $(\Ai)_{n}$ up to $n=5$ and we find the same result as the Saneblidze-Umble diagonal. So it is reasonable to conjecture that they coincide.

In the last part we give a similar interpretation of the map $p:\Omega B As{\to}  \Ai$ constructed in \cite{MS} and giving rise to the Saneblidze-Umble diagonal. It is induced by the deformation of the ``main cube'' into the big cell.

Moreover we provide an explicit comparison map between the two  models $\Omega B As$ and $\AAi $ by using the simplicialization of the cubical decomposition of the Stasheff polytope. It should prove useful in the comparison of the two diagonals.
\vskip 1cm

\noindent {\bf Acknowledgement} I thank Bruno Vallette for illuminating discussions on  the algebras up to homotopy and Samson Saneblidze for sharing his drawings with me some years ago. Thanks to Emily Burgunder, Martin Markl, Samson Saneblidze, Jim Stasheff and Ron Umble for their comments on the previous versions of this paper.

This work is partially supported by the French agency ANR.

\section{Stasheff polytope (associahedron)} We recall briefly the construction of the Stasheff polytope,  also called associahedron, and its simplicial realization, which is the key tool of this paper.
All chain complexes in this paper are made of free modules over a commutative ring $\KK$ (which can be $\ZZ$ or a field).

\subsection{Planar binary trees}\label{planarbinarytree}  We denote by $PBT_{n}$ the set of \emph{planar binary trees} having $n$ leaves:
$$PBT_{1}:=\{ | \}, PBT_{2}:=\{\arbreA \}, PBT_{3}:=\{\arbreAB , \arbreBA \}, $$
\begin{displaymath}
PBT_{4}:=\{\arbreABC, \arbreBAC, \arbreACA, \arbreCAB, \arbreCBA \}\ .
\end{displaymath}

So $t\in PBT_{n}$ has one root, $n$ leaves, $(n-1)$ internal vertices, $(n-2)$ internal edges. Each vertex is binary (two inputs, one output). The number of elements in $PBT_{n+1}$ is known to be the \emph{Catalan number}\index{Catalan number} $c_{n}= \frac{(2n)!}{n!\, (n+1)!}$. There is a partial order on $PBT_{n}$, called the \emph{Tamari order}, defined as follows. On $PBT_{3}$ it is given by 
$$\arbreAB \to \arbreBA\ .$$
More generally, if $t$ and $s$ are two planar binary trees with the same number of leaves, there is a covering relation $t \to s$ if and only if $s$ can be obtained from $t$ by replacing a local pattern like $\arbreAB$ by $\arbreBA$. In other words $s$ is obtained from $t$ by moving a leaf or an internal edge from left to right over a fork.

Examples:
$$\KzeroTree \qquad \KunTree \qquad \KdeuxTree$$
where the elements of $PBT_{4}$ (listed above) are denoted $123, 213, 141, 312, 321,$ respectively. We recall from \cite{JLLStasheff} how this way of indexing is obtained.  First we label the leaves of a tree from left to right by $0,1,2,\ldots $ . Then we label the vertices by $1,2,\ldots $ by saying that the label $i$ vertex lies in between the leaves $i-1$ and $i$ (drop a ball). To any binary tree $t$ we associate a sequence of integers $x_{1}x_{2}\cdots x_{n-1}$ as follows: $x_{i}= a_{i}b_{i}$ where $a_{i}$ (resp. $b_{i}$) is  the number of leaves on the left (resp. right) side of the $i$th vertex.

\subsection{Shortest path and long path}
The Tamari poset admits an initial element: the left comb $12 \cdots (n-1)$, and a terminal element: the right comb $(n-1) (n-2)\cdots 1$. There is a \emph{shortest path} from the initial element to the terminal element. It is made of the trees which are the grafting of some left comb with a right comb. In $PBT_{n}$ there are $n-1$ of them. This sequence of planar binary trees will play a significant role in the comparison of different cell realizations of the Stasheff polytope.

Example: the shortest path in $PBT_{4}$:
$$\arbreABC \to \arbreACA \to \arbreCBA$$

We also define ``the long path'' as follows. The \emph{long path} from the left comb to the right comb is obtained by taking a covering relation at each step with the following rule: the vertex which is moved is the one with the smallest label (among the movable vertices, of course).

Examples: $n=2$

$$\arbreABC \to \arbreBAC \to \arbreCAB \to \arbreCBA$$
$n=3$

$\begin{array}{l c}
 \arbreABCD \to \arbreBACD \to \arbreCABD \to \arbreDABC  \to &\qquad
\end{array}$

$\begin{array}{c r}
\qquad & \arbreDBAC \to \arbreDCAB \to \arbreDCBA
\end{array}$

Observe that there are (for $n\geq 3$) other paths with the same length.

\subsection{Planar trees}\label{planartree}   We now consider the planar trees for which an internal vertex has one root and $k$ leaves, where $k$ can be any integer greater than or equal to 2. We denote by $PT_{n}$ the set of planar trees with $n$ leaves:
$$PT_{1}:=\{ | \}, PT_{2}:=\{\arbreA \}, PT_{3}:=\{\arbreAB , \arbreBA, \arbreBB \}, $$
$$PT_{4}:=\{\arbreABC, \ldots , \arbreBBC, \ldots, \arbreCCC \}\ .$$
  Each set $PT_{n}$ is graded according to the number of internal vertices, i.e.~$PT_{n}= \bigcup_{p=1}^{p=n} PT_{n,p}$ where $PT_{n,p}$ is the set of planar trees with $n$ leaves and $p$ internal vertices. For instance $PT_{n,1}$ contains only one element which we call the \emph{$n$-corolla}\index{corolla} (the last element in the above sets). It is clear that $PT_{n,n-1}= PBT_{n}$. 
     
    We order the vertices of a planar tree by using the same procedure as for the planar binary trees.

  \subsection{The Stasheff polytope,  alias associahedron}\label{associahedron}  The \emph{associahedron} is a cellular complex $\KKK^{n}$ of dimension $n$, first constructed by Jim Stasheff \cite{Stasheff}, which can be realized as a convex polytope whose extremal vertices are in one-to-one correspondence with the planar binary trees in $PBT_{n+2}$. We showed in \cite{JLLStasheff} that it is the convex hull of the points $M(t)=(x_{1}, \ldots , x_{n+1})\in \RR^{n+1}$. The edges of the polytope are indexed by the covering relations of the Tamari poset. 
  
  Examples:
  
  $$\KzeroF \qquad \KunF\qquad\quad \KdeuxF\qquad\quad \KtroisF
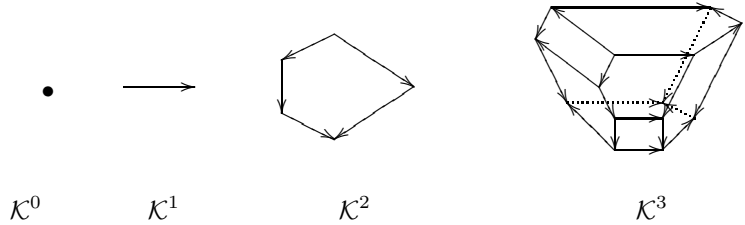
$$
  $$\KKK^0 \qquad \qquad  \KKK^1 \qquad  \qquad  \qquad \KKK^2 \qquad \qquad \qquad  \qquad  \qquad \KKK^3\qquad \qquad$$
  
   Its $k$-cells are in one-to-one correspondence with the planar trees in $PT_{n+2, n+1-k}$. For instance the $0$-cells are indexed by the planar binary trees, and the top cell is indexed by the corolla.
   
It will prove helpful to adopt the notation $\KKK^t$ to denote the cell in $\KKK^n$ indexed by $t\in PT_{n+2}$. For instance, if $t$ is the corolla, then $\KKK^t= \KKK^n$. As a space $\KKK^t$ is the product of $p$ associahedrons (or associahedra, as you like), where $p$ is the number of internal vertices of $t$: 
  $$\KKK^t = \KKK^{i_{1}}\times \cdots \times \KKK^{i_{p}}$$
  where $i_{j}+2$ is the number of inputs of the $j$th internal vertex of $t$.  For instance, if $t= \arbreBBDD$, then $\KKK^t= \KKK^1\times \KKK^1$.

The shortest path and the long path defined combinatorially in \ref{planarbinarytree} give rise to concrete paths on the associahedron.

To the cellular complex $\KKK^n$ we associate its chain complex $C_{*}(\KKK^n)$. The module of $k$-chains admits the set of trees $PT_{n+2,n+1-k}$ as a basis:
$$C_{k}(\KKK^n)= \KK[PT_{n+2,n+1-k}].$$
In particular $C_{0}(\KKK^n)= \KK[PBT_{n+2}]$ and $C_{n}(\KKK^n)= \KK\ t_{n+2}$ where $t_{n+2}$ is the corolla.

\subsection{The simplicial associahedron}\label{simplicialSP} In \cite{JLLparking} we constructed a simplicial set $\Ksimp n $ whose geometric realization gives a simplicial decomposition of the associahedron. In other words the associahedron $\KKK^n$ is viewed as a union of $n$-simplices (there are $(n+1)^{n-1}$ of them). 
This simplicial decomposition is constructed inductively as follows. We fatten the simplicial set $\Ksimp {n-1}$ into a new simplicial set $fat\Ksimp {n-1}$, cf. \cite{JLLparking}. Then $\Ksimp n$ is defined as the cone over $fat\Ksimp {n-1}$ (as in the original construction of Stasheff \cite{Stasheff}). 

For $n=1$, we have $\Ksimp 1 = \KKK^1 =[0,1]$ (the interval).

Examples:  $\Ksimp 2$ and $fat\Ksimp 3$
$$\KdeuxTfat \qquad \KtroisTfat$$
 Since, in the process of fattenization, the new cells are products of smaller dimensional associahedrons we get the following main property.

\begin{proposition}\label{mainproperty} The simplicial decomposition of a face  $\KKK^{i_{1}}\times \cdots \times \KKK^{i_{k}}$ of $\KKK^{n}$ is the product of the simplicializations of each component $\KKK^{i_{j}}$.
\end{proposition}
\begin{proo} It is immediate from the inductive procedure which constructs $\KKK^{n}$ out of $\KKK^{n-1}$.
\end{proo}

Considered as a cellular complex, still denoted $\Ksimp n $, the simplicialized associahedron gives rise to a chain complex denoted $C_{*}(\Ksimp n )$. This chain complex is the normalized chain complex of the simplicial set. It is the quotient of the chain complex associated to the simplicial set, divided out by the degenerate simplices (cf.~for instance \cite{MacLane} Chapter VIII). A basis of $C_{0}(\Ksimp n )$ is given by $PBT_{n+2}$ and a basis of $C_{n}(\Ksimp n )$ is given by the $(n+1)^{n-1}$ top simplices (in bijection with the parking functions, cf.~\cite{JLLparking}). It is zero higher up.

In the sequel ``a simplex of $\Ksimp n$'' always mean a nondegenerate simplex of $\Ksimp n$.

Among the top simplices there is a particular one which we call the \emph{main simplex}. Its vertices are indexed by the planar binary trees which are part of the shortest path constructed in \ref{planarbinarytree} (observe that the shortest path has $n+1$ vertices).

Examples (the main simplex is highlighted):

$$\KdeuxMainS \qquad \KtroisMainS$$

\section{The operad $\AAi$}\label{AA-infinity} We construct the operad $\AAi$ and we construct a diagonal on it. A morphism from the operad $\Ai$ governing the associative algebras up to homotopy to the operad $\AAi$ is deduced from the simplicial structure of the associahedron. 

\subsection{Differential graded non-symmetric operad \cite{MSS}} By definition a \emph{ differential graded non-symmetric operad}, dgns operad for short, is a family of chain complexes  $\PP_{n} = (\PP_{n},d)$ equipped with chain complex morphisms 
$$\cc_{i_{1}\cdots i_{n}}: \PP_{n}\t \PP_{i_{1}}\t \cdots\t \PP_{i_{n}} \to \PP_{i_{1}+\cdots +i_{n}},$$
which satisfy the following associativity property. Let $\PP$ be the endofunctor of the category of chain complexes over $\KK$ defined by $\PP(V) := \bigoplus_{n} \PP_{n} \t V^{\t n}$. The maps 
$\cc_{i_{1}\cdots i_{n}}$ give rise to a transformation of functors $\cc : \PP \circ \PP \to \PP$. This transformation of functors $\cc$ is supposed to be associative. Moreover we suppose that $\PP_{0}= 0, \PP_{1} = \KK$ (trivial chain complex concentrated in degree 0). The transformation of functors $\Id \to \PP$ determined by $\PP_{1}$ is supposed to be a unit for $\cc$. So we can denote by $\id$ the generator of $\PP_{1}$. Since $\PP_{n}$ is a graded module, $\PP$ is bigraded. The integer $n$ is called the ``arity'' in order to differentiate it from the degree of the chain complex.

\subsection{The fundamental example $\Ai$} The operad $\Ai$ is a dgns operad constructed as follows:
$$A_{\infty , n}:= C_{*}(\KKK^{n-2})\ ({\rm chain\ complex\ of\ the\ cellular\ space\ } \KKK^{n-2}).$$

Let us denote by $As\ac$ the family of one dimensional modules $(As\ac_{n})_{n\geq 1}$ generated by the corollas (unique top cells). It is easy to check that there is a natural identification of graded (by arity) modules $\Ai = \TTT(As\ac)$, where $\TTT(As\ac)$ is the free ns operad over $As\ac$. This identification is given by grafting on the leaves as follows. Given trees $t, t_{1}, \ldots , t_{n}$ where $t$ has $n$ leaves, the tree $\cc(t; t_{1}, \ldots , t_{n})$ is obtained by identifying the $i$th leaf of $t$ with the root of $t_{i}$. For instance:
$$\cc(\arbreA ; \arbreBA , \arbreA) = \arbreBADA\ .$$
Moreover, under this identification, the composition map $\cc$ is a chain map, therefore $\Ai$ is a dgns operad.

This construction is a particular example of the so-called ``cobar construction'' $\Omega$, i.e.~$\Ai = \Omega As\ac $ where $As\ac$ is considered as the cooperad governing the coassociative coalgebras (cf.~\cite{MSS}).

For any chain complex $A$ there is a well-defined dgns operad $\End(A)$ given by $\End(A)_{n}= \Hom(A^{\t n} , A)$. An $\Ai$-algebra is nothing but a morphism of operads $\Ai \to \End(A)$. The image of the corolla under this isomorphism is the $n$-ary operation $\mu_{n}$ alluded to in the introduction.

\subsection{Hadamard product of operads, the diagonal problem} Given two operads $\PP$ and $\QQQ$, their Hadamard product, also called tensor product, is the operad $\PP \t \QQQ$ defined as $(\PP \t \QQQ)_{n} := \PP_{n} \t \QQQ_{n}$. The composition map is simply the tensor product of the two composition maps.

It is a long-standing problem to decide if, given two $\Ai$-algebras $A$ and $B$, there is a natural $\Ai$-structure on their tensor product $A\t B$ which extends the natural dg nonassociative algebra structure, cf.~\cite{Proute, GZ}. It amounts to construct a diagonal on $\Ai$, i.e.~an operad morphism $\DD: \Ai \to \Ai\t \Ai$, since, by composition, we get an $\Ai$-structure on $A\t B$:
$$\Ai \to \Ai\t \Ai\to \End(A)\t \End(B) \to \End(A\t B)\ .$$
Let us recall that the classical associative structure on the tensor product of two associative algebras can be interpreted operadically as follows. There is a diagonal on the operad $As$ given by 
$$As_{n} \to As_{n} \t As_{n}, \quad \mu_{n} \mapsto \mu_{n} \t \mu_{n} \ .$$
Since we want the diagonal $\DD$  to be compatible with the diagonal on $As$, there is no choice in arity 2, and we have $\DD(\mu_{2}) = \mu_{2}\t \mu_{2}$.  Observe that these two elements are in degree $0$. In arity $3$, since $\mu_{3}$ is of degree $1$ and $\mu_{3}\t \mu_{3}$ of degree $2$, this last element cannot be the answer. In fact there is already a choice (parameter $a$) for a solution:
\begin{eqnarray*}
\DD(\arbreBB ) &=&a\Big(\arbreBB \t \arbreAB +\arbreBA \t \arbreBB \Big)\\
& & + (1-a)\Big(\arbreBB \t \arbreBA +\arbreAB \t \arbreBB \Big) .\\
\end{eqnarray*}

By some tour de force Samson Saneblidze and Ron Umble constructed such a diagonal on $\Ai$  in \cite{SU1}. Their construction was re-interpreted in \cite{MS} by Markl and Shnider through the Boardman-Vogt construction (see section \ref{MS} below for a brief account of their work). We will use the simplicialization of the associahedron described in \cite{JLLparking} to give the solution to the diagonal problem.

\subsection{Construction of the operad $\AAi$}\label{AAinfini} We define the dgns operad $\AAi$ as follows. The chain complex $\Ainf n$ is the chain complex of the simplicialization of the associahedron considered as a cellular complex (cf.~\ref{simplicialSP}):
$$\Ainf n = C_{*}(\Ksimp {n-2})\ .$$
In low dimension we take $\Ainf 0 = 0, \Ainf 1= \KK\, \id$. So a basis of $\Ainf n $ is made of the (nondegenerate) simplices of $\Ksimp {n-2}$.
Let us now construct the composition map 
$$\cc=\cc^{\AAi}: \Ainf n \t \Ainf {i_{1}} \t   \cdots \t \Ainf {i_{n}} \to \Ainf {i_{1}+\cdots +i_{n}}.$$

We denote by $\Delta^{k}$ the standard $k$-simplex. Let $\iota: \Delta^{k}\mono \Ksimp {n-2}$ be a cell, i.e.~a linear generator of $C_{k}(\Ksimp {n-2})$. Given such cells
$$\iota_{0}\in  \Ainf n ,\ \iota_{1}\in  \Ainf {i_{1}} ,\ \ldots ,\   \iota_{n}\in  \Ainf {i_{n}}$$
we construct their image $\cc(\iota_{0}; \iota_{1} , \ldots ,  \iota_{n})\in  \Ainf {m}$, where $m:= i_{1}+\cdots +i_{n}$ as follows. We denote by $k_{i}$ the dimension of the cell $\iota_{i}$.

Let $t_{n}$ be the $n$-corolla in $PT_{n}$ and let $s:= \cc(t_{n}; t_{i_{1}} , \ldots ,  t_{i_{n}})\in PT_{m}$ be the grafting of the trees $t_{i_{1}} , \ldots ,  t_{i_{n}}$ on the leaves of $t_{n}$. As noted before this is the composition in the operad $\Ai$. The tree $s$ indexes a cell $\KKK^s$ of the space $\KKK^{m-2}$, which is combinatorially homeomorphic to $\KKK^{n-2}\times \KKK^{i_{1}-2}\times \cdots \times \KKK^{i_{n}-2}$. In other words it determines a map 
$$s_{* }:\KKK^{n-2}\times \KKK^{i_{1}-2}\times \cdots \times \KKK^{i_{n}-2}= \KKK^s \mono  \KKK^{m-2}.$$
The product of the inclusions $\iota_{j}, j=0,\ldots , n,$ defines a map
$$\iota_{0}\times \iota_{1}\times \cdots \times \iota_{n}: \DD^{k_{0}}\times \DD^{k_{1}}\times
\cdots \times  \DD^{k_{n}} \mono  \KKK^{n-2}\times \KKK^{i_{1}-2}\times \cdots \times \KKK^{i_{n}-2}.$$
Let us recall that a product of standard simplices can be decomposed into the union of standard simplices. These pieces are indexed by the multi-shuffles $\aa$. Example: $\DD^1 \times \DD^1 = \DD^2 \cup \DD^2$:
$$
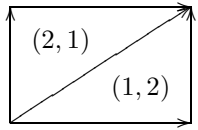
\xymatrix@R=2pt@C=2pt{
*{}\ar[rrr] & & & *{}\\
 & (2,1) & & \\
  & & (1,2) & \\
*{}\ar[uuu]\ar[rrr] \ar[uuurrr] & & & *{}\ar[uuu] \\
}$$
So, for any multi-shuffle $\aa$ there is a map
$$f_{\aa}: \DD^{l} \to  \DD^{k_{0}}\times \DD^{k_{1}}\times
\cdots \times  \DD^{k_{n}},$$
where $l= k_{0}+\cdots + k_{n}$. 
By composition of maps we get
$$s_{*} \circ (\iota_{0}\times \cdots \times \iota_{n})\circ f_{\aa}:  \DD^{l} \to  \KKK^{m-2}$$
which is a linear generator of $C_{l}(\Ksimp {m-2})$ by construction of the triangulation of the associahedron, cf.~\cite{JLLStasheff}. By definition 
$\cc(\iota_{0}; \iota_{1} , \ldots ,  \iota_{n})$ is the algebraic sum of the cells $s_{*} \circ (\iota_{0}\times \cdots \times \iota_{n})\circ f_{\aa}$ over the multi-shuffles.

\begin{proposition} The graded chain complex $\AAi$ and $\cc$ constructed above define a dgns operad, denoted $\AAi$. The operad $\AAi$ is a model of the operad $As$.
\end{proposition}
\begin{proo} We need to prove associativity for $\cc$. It is an immediate consequence of the associativity for the composition of trees (operadic structure of $\Ai$) and the associativity property for the decomposition of the product of simplices into simplices.

Since the associahedron is contractible, taking the homology gives a graded linear map 
$C_{*}(\Ksimp {n-2})\to \KK\, \mu_{n}$, where $\mu_{n}$ is in degree 0. This map obviously induces an isomorphism on homology. These maps assemble into a dgns operad morphism $\AAi \to As$. Since it  is a quasi-isomorphism, $\AAi$ is a resolution of $As$, that is a model of $As$ in the category of dgns operads.
\end{proo}


\begin{proposition} The operad $\AAi$ admits a coassociative diagonal.
\end{proposition}

\begin{proo} This diagonal $\DD : \AAi \to \AAi \t \AAi$ is determined by its value in arity $n$ for all $n$, that is a chain complex morphism
$$C_{*}(\Ksimp {n-2})\to C_{*}(\Ksimp {n-2} )\t C_{*}(\Ksimp {n-2} ).$$
This morphism is defined as the composite
$$C_{*}(\Ksimp {n-2})\stackrel{\DD_{*}}{\longrightarrow} C_{*}(\Ksimp {n-2} \times \Ksimp {n-2}) \stackrel{AW}{\longrightarrow} C_{*}(\Ksimp {n-2})\t C_{*}(\Ksimp {n-2}),$$
where $\DD_{*}$ is induced by the diagonal on the simplicial set, and where $AW$ is the Alexander-Whitney map. Let us recall from \cite{MacLane}, Chapter VIII, the construction of the AW map. Denote by $d_{0}, \ldots , d_{k}$ the face operators of the simplicial set. If $x$ is a simplex of dimension $k$, then we define $d_{max}(x):= d_{k}(x)$. So, for instance $(d_{max})^2(x) = d_{k-1}d_{k} (x)$. By definition the AW map on $C_{k}$ is given by
$$(x,y) \mapsto \sum_{i=0}^k \big( (d_{max})^{k-i} (x), (d_{0})^{i}(y) \big).$$

It is straightforward to check that this diagonal is compatible with the operad structure.

The coassociativity property follows from the coassociativity property of the Alexander-Whitney map.
\end{proo}

\subsection{Comparing $\Ai$ to $\AAi$} Since $\Ksimp {n}$ is a decomposition of $\KKK^{n}$, there is a chain complex map
$$q': C_{*}(\KKK^{n}) \to C_{*}(\Ksimp {n}),$$
where a cell of $\KKK^{n}$ is sent to the algebraic sum of the simplices it is made of. 

\begin{proposition}  The map $q':\Ai \to \AAi$ induced by the maps $q': C_{*}(\KKK^{n}) \to C_{*}(\Ksimp {n})$ is a quasi-isomorphism of dgns operads.
\end{proposition}
\begin{proo} It is sufficient to prove that the maps $q'$ on the chain complexes are compatible with the operadic composition:
$$q'(\cc^{As}(t; t_{1},\ldots , t_{n}))= \cc^{\AAi}(q'(t); q'(t_{1}),\ldots , q'(t_{n})).$$
This equality follows from the definition of $\cc^{\AAi}$ given in \ref{AAinfini} and Proposition \ref{mainproperty}.
\end{proo}

 Moreover we have commutative diagrams:

$\xymatrix@R=20pt@C=20pt{
C_{*}(\KKK^{n-2}) \ar[rr]^-{q'}\ar[dr]^-{H_{*}}  && C_{*}(\Ksimp {n-2})\ar[dl]_-{H_{*}}\\
&\KK\, \mu_{n}& 
}$\qquad 
$\xymatrix@R=20pt@C=20pt{
\Ai \ar[rr]^-{q'}\ar[dr]^-{H_{*}}   && \AAi\ar[dl]_-{H_{*}}\\
&As& 
}$

\section{From $\AAi$ to $\Ai$} 
The aim of this section is to construct a quasi-inverse to $q'$, that is a quasi-isomorphism of dgns operads $p': \AAi\to \Ai$. We first construct chain maps  $p': C_{*}(\Ksimp {n}) \to C_{*}(\KKK^{n})$ by using a deformation of the main simplex to the top cell of the associahedron.

\subsection{Deformation of the cube}\label{deformedCube} The cube $I^n$ is a polytope whose vertices are indexed by $(x_{1}, \ldots, x_{n})$, where $x_{i}= 0$ or $1$.  The long path in $I^n$ is, by definition,  the path 
$$(0,\ldots , 0,0) \to (0,\ldots , 0,1) \to (0,\ldots , 1,1) \to \cdots \to (1,\ldots , 1,1) .$$
  The cube is a cell complex which can be decomposed into $n!$ top simplices, i.e.~viewed as the realization of a simplicial set $I^n_\textrm{simp}$. The simplex which corresponds to the identity permutation is called the main simplex of the cube. Let us describe the deformation from the main simplex to the cube, which gives rise to a chain map
$$p': C_{*}(I^n_\textrm{simp}) \to C_{*}(I^n)\ .$$
We work by induction on $n$. In $I^2_\textrm{simp}$ the main simplex, denoted $\aa$,  is deformed to the square by pushing the diagonal to the long path:

$$
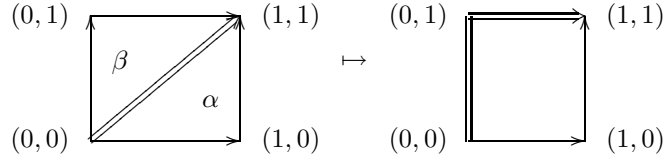
\xymatrix@R=2pt@C=2pt{
(0,1)&*{}\ar[rrr] & & & *{}&(1,1)                                            &                          &(0,1)&  *{}\ar@{=>}[rrr] & & & *{}&(1,1)\\
& & \bb\quad & &                   &                                                          & \mapsto &  &                & &\qquad  &    &\\
&   & &\quad \aa &          &                                                  &                              &&&                 & &   &\\
(0,0)&*{}\ar[uuu]\ar[rrr] \ar@{=>}[uuurrr] & & & *{}\ar[uuu]&(1,0)  &                 & (0,0)&*{}\ar@{=}[uuu]\ar[rrr]  & &&*{}\ar[uuu]&(1,0)
}$$

So $p'$ is given by the identity on the boundary and by

$$\begin{array}{c c c}
((0,0),(1,1))& \mapsto& ((0,0),(0,1)) + ((0,1),(1,1))  \\
\aa&\mapsto& I^2  \\
\bb&\mapsto& 0  \\
\end{array}$$
on the interior simplices. So, under this deformation, the main simplex $\aa$ is mapped to the whole square and the other simplex $\bb$ is flattened. More generally, the main simplex of $I^n_\textrm{simp}$ is deformed into the top cell of $I^n$ by sending the diagonal to the long path. The other edges of the main simplex are deformed according to the lower dimensional deformation.

\subsection{Deformation of a product of simplices}\label{deformedProduct}  Similarly we define a deformation of the product of simplices $\DD^r\times \DD^s$ as follows. Let us denote by $\{\u{0},\ldots , \u{r}\}$ the vertices of $\DD^r$. The \emph{main simplex} of $\DD^r\times \DD^s$ is chosen as being the simplex $\DD^{r+s}$ with vertices 
$$(\u{0},\u{0}), (\u{1},\u{0}),\ldots , (\u{r},\u{0}),(\u{r},\u{1}),\ldots , (\u{r},\u{s}),\ .$$
We deform the main simplex into the whole product by induction on $s$. So it suffices to give the image of the edge $((\u{0},\u{0}), (\u{r},\u{s}))$. We send it to the ``long path'' defined as
$$(\u{0},\u{0}), (\u{0},\u{1}),\ldots , (\u{0},\u{s}),(\u{1},\u{s}),\ldots , (\u{r},\u{s}),\ .$$
Under this deformation the main simplex becomes the whole product and all the other simplices are flattened. This deformation defines a chain complex morphism
$$C_{*}((\DD^r\times \DD^s)_{\textrm{simp}}) \to C_{*}(\DD^r\times \DD^s)$$
where, on the right side, $\DD^r\times \DD^s$ is considered as a cell complex with only one $r+s$ cell.

Observe that some simplices may happen to be deformed into cells of various dimensions. For instance in $\DD^2\times \DD^1$ the triangle with vertices $(\u{0},\u{0}), (\u{1},\u{1}),  (\u{2},\u{1})$ is deformed into the union of the edge $(\u{0},\u{0})\to  (\u{0},\u{1})$ and the triangle with vertices $(\u{0},\u{1}), (\u{1},\u{1}),  (\u{2},\u{1})$.

\subsection{Deformation of the associahedron}\label{deformedAssociahedron} We construct a topological deformation of the simplicial associahedron by pushing the main simplex to the whole associahedron. All the other simplices are going to be flattened. This topological deformation will induce the chain map $p'$ we are looking for. The process is analogous to what we did for the cube and the product of simplices above. We work by induction on the dimension.

For $n=1$, there is no deformation since $\Ksimp 1 = \KKK^1$. For $n=2$ the deformation is the identity on the boundary and the only edge of the main simplex which is not on the boundary is ``pushed'' to the long path.  

$$\KdeuxSimp \qquad   \mapsto\qquad  \KdeuxPushed$$

In the meantime, the other interior edge is pushed to the union of two boundary edges and the two other top simplices are flattened.

For higher $n$ we use the inductive construction of $\Ksimp n$ out of $fat\Ksimp {n-1}$. We suppose that the deformation is known for any $i<n$ and we construct it on $fat\Ksimp {n-1}$. The simplicial set $fat\Ksimp {n-1}$ is the union of the simplicial sets of the form $\Ksimp t = \Ksimp i \times \Ksimp j$ indexed by some trees $t$ with one and only one internal edge. The main simplex of this product is the main simplex $\DD^{i+j}$ of $\DD^{i} \times \DD^j$ where $\DD^{i}$, resp.~$\DD^{j}$, is the main simplex of $\Ksimp i $, resp.~ $\Ksimp j $. The deformation is obtained by, first, deforming  the main simplex $\DD^{i+j}$ into $\DD^{i} \times \DD^j$ as described in \ref{deformedProduct} and then use the inductive hypothesis (deformation from the main simplex to the associahedron).

The deformation of the interior cells is obtained by pushing the main simplex of $\Ksimp n$ to the top cell. It is determined by the image of the edges of the main simplex. By induction, it suffices to construct the image of the edge which goes from the vertex indexed by the left comb (initial element) to  the vertex indexed by the  right comb (terminal element). We choose to deform it to the long path of the associahedron. Since any simplex of $\Ksimp n$ is either on the boundary, or is a cone (for the last vertex) over a simplex in the boundary, we are done.
In particular, the edge going from a 0-simplex labelled by the tree $t$ to the right comb is deformed into a path made of 1-cells of the associahedron, constructed with the same rule as in the construction of the long path.

The deformed tetrahedron:

$$\KtroisDeformed$$

\subsection{The map $p': C_{*}(\Ksimp {n}) \to C_{*}(\KKK^{n})$} 

We define the map $p'$ as follows. 
Under the deformation map any simplex of $\Ksimp {n}$ is sent to the union of cells of $\KKK^{n}$. The image of such a simplex  under $p'$  is the algebraic sum of the cells of the same dimension in the union. For instance, the main simplex is sent to the top cell (indexed by the corolla), and all the other top simplices are sent to $0$, since under the deformation they are flattened. From its topological nature it follows that $p'$ is a chain complex morphism.

In low dimension we get the following.
For $n=1$, the map $p'$ is the identity. For $n=2$, the map $p'$ is the identity on the $0$-simplices and the $1$-simplices of the boundary, and on the interior cells, we get:
$$\begin{array}{c c c}
a=\Big(\arbreABC, \arbreACA, \arbreCBA\Big)& \mapsto& \arbreCCC \\
b=\Big(\arbreABC, \arbreBAC, \arbreCBA\Big)& \mapsto& 0 \\
c=\Big(\arbreBAC, \arbreCAB, \arbreCBA\Big)& \mapsto& 0 \\
\end{array}$$
$$\begin{array}{c c c}
\Big(\arbreABC,  \arbreCBA\Big)& \mapsto& \arbreBBC + \arbreCAC +\arbreCBB \\
\Big(\arbreBAC,  \arbreCBA\Big)& \mapsto& \arbreCAC +\arbreCBB \\
\end{array}$$

 Here are examples of the image under $p'$ of  some interior 2-dimensional simplices for $n=3$:
$$\Big(\arbreABCD, \arbreADBA, \arbreDCBA\Big) \mapsto\qquad\qquad$$
$$\qquad\qquad \arbreCCCD + \arbreDDAD + \arbreDDBB .$$
$$\Big(\arbreBACD, \arbreBADA, \arbreDCBA\Big) \mapsto\qquad\qquad$$
$$\qquad\qquad \arbreDADD + \arbreDBBD + \arbreDDAD .$$

\begin{proposition} The chain maps $p': C_{*}(\Ksimp {n}) \to C_{*}(\KKK^{n})$ assemble into a morphism of dgns operads $p':\AAi \to \Ai$.
\end{proposition}
\begin{proo} 
We adopt the notation of \ref{AAinfini} where the operadic composition map $\cc^{\AAi}$ is constructed. From this construction it follows that there is a main simplex in $\oo:= \cc^{\AAi}(\iota_{0}; \iota_{1} , \ldots ,  \iota_{n})$ if and only if all the simplices $\iota_{j}$ are main simplices. 

Supppose that one of them, say $\iota_{j}$ is not a main simplex. Then we have $p'(\iota_{j})=0$, and therefore 
 $\cc^{\Ai}(p'(\iota_{0}); p'(\iota_{1}) , \ldots ,  p'(\iota_{n}))=0$. But since there is no main simplex in $\oo$, we also get $p'(\oo)=0$ as expected.
 
 Suppose that all the simplices are main simplices. Then $p'(\iota_{j})= t_{i_{j}}$ for all $j$ and therefore
 $\cc^{\Ai}(p'(\iota_{0}); p'(\iota_{1}) , \ldots ,  p'(\iota_{n}))=\cc^{\Ai}(t_{i_{0}}; t_{i_{1}}, \ldots , t_{i_{n}})$. On the other hand 
 $\oo$ contains the main simplex, therefore $p'(\oo)= \cc^{\Ai}(t_{i_{0}}; t_{i_{1}}, \ldots , t_{i_{n}})$ and we are done.
\end{proo}

\begin{corollary} The composite
$$ \Ai \stackrel{q'}{\to} \AAi \stackrel{\DD}{\longrightarrow}  \AAi \t \AAi \stackrel{p'\t p'}{\longrightarrow} \Ai \t \Ai\ .$$
is a diagonal for the operad $\Ai$.
\end{corollary} 
\begin{proo} It is immediate to check that this composite of dgns operad morphisms sends $\mu_{2}$ to $\mu_{2}\t \mu_{2}$, since $\mu_{2}$ corresponds to the 0-cell of $\KKK^0$.
\end{proo}

\begin{proposition} If $A$ is an associative algebra and $B$ an $\Ai$ algebra, then the $\Ai$-structure on $A\t B$ is given by 
$$\mu_{n}(a_{1}\t b_{1}, \ldots , a_{n}\t b_{n}) = a_{1}\cdots a_{n}\t \mu_{n}(b_{1}, \ldots , b_{n}) .$$
\end{proposition}
\begin{proo} In the formula for $\DD$ we have $\mu_{n}=0$ for all $n\geq 3$, that is, any tree with a $k$-valent vertex for $k\geq 3$ is $0$ on the left side. Hence the only term which is left is $comb \t corolla$, whence the assertion.
\end{proo}

\subsection{The first formulas} Let us give the explicit form of $\DD(\mu_{n})$ for $n=2,3,4$:
\begin{eqnarray*}
\DD(\arbreA) &=& \arbreA \t \arbreA ,\\
\DD\big(\arbreBB\big) &=& \arbreAB \t \arbreBB +  \arbreBB \t \arbreBA ,
\end{eqnarray*}
$$\DD\Big(\arbreCCC\Big) = \arbreABC \!\t\! \arbreCCC +  \arbreCCC \!\t\! \arbreCBA $$
$$+\!\!\arbreACC\!\t\! \arbreCCA\!\! -\!\!\arbreBBC \!\t\! \arbreCBB\!\! - \!\! \arbreBBC \!\t\! \arbreCAC\!\! -\!\!\arbreCAC\! \t\! \arbreCBB\! .$$

In this last formula the first three summands comes from the triangle $(123, 141, 321)$, the next two summands come from the triangle $(123, 213, 321)$ and the last summand comes from the last triangle $(213, 312, 321)$. It is exactly the same formula as the one obtained by Saneblidze and Umble (cf.~\cite{SU1} example 1, \cite{MS} exercise 12).
Topologically the diagonal of the pentagon is approximated as a union of products of cells as follows:

$$\KdeuxDIAG$$
Each cell of this decomposition corresponds to a summand of the above formula, which indicates where the cell goes in the product $\KKK^2\times \KKK^2$.

\subsection{On the non-coassociativity of the diagonal} Though the diagonal of $\AAi$ that we constructed is coassociative, the diagonal of $\Ai$ is not. In fact it has been shown in \cite{MS} that there does not exist any coassociative diagonal on $\Ai$. The obstruction to coassociativity can be seen topologically on the picture ``Iterated diagonal''.

\begin{figure}[htbp]
\begin{center}
\input{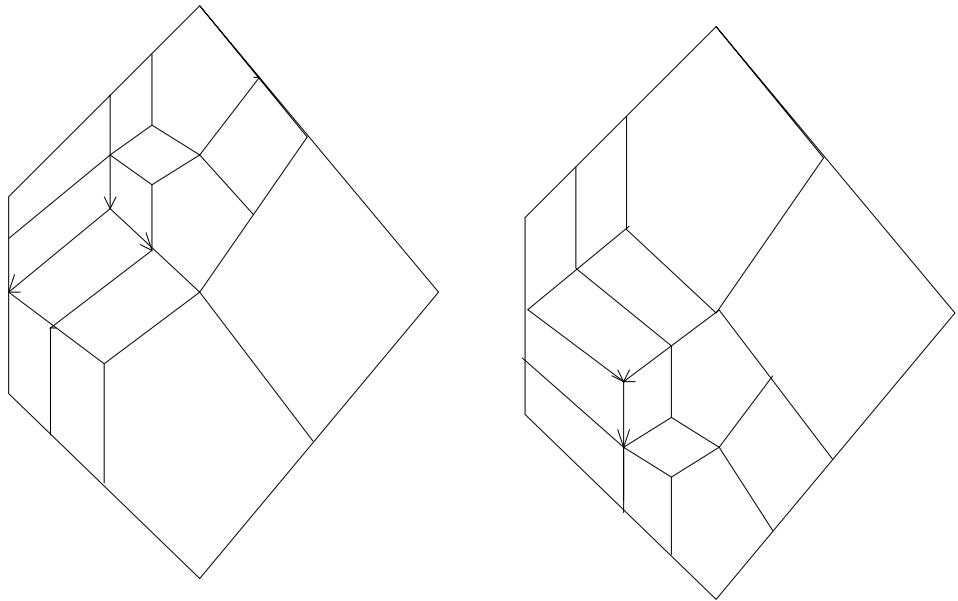}
\caption{Iterated diagonal}
\label{figure:assoc}\end{center}
\end{figure}

Both pictures are the same combinatorially, except for an hexagon (highlighted on the pictures), which is the union of 3 squares one way on the left and the other way on the right. This is the obstruction to coassociativity. Of course there is a way to reconcile these two decompositions via a homotopy which is given by the cube.

Exercise 1. Show that the image of this cube in $\KKK^2\times \KKK^2\times \KKK^2$ is indexed by $\arbreBBC \times \arbreCAC\times  \arbreCBB$.

Exercise 2. Compare the five iterated diagonals of the next step (some nice pictures to draw).

\section{Comparing the operads $\AAi$ and $\Omega B As$}\label{MS} We first give a brief account of \cite{MS, SU1} where a diagonal of the operad $\Ai$ is constructed by using a coassociative diagonal on the dgns operad $\Omega B As$. Then we compare the two operads $\AAi$ and $\Omega B As$.

\subsection{Cubical decomposition of the associahedron \cite{BV}} The associahedron can be decomposed into cubes as follows. 

For each tree $t\in PBT_{n+2}$ we take a copy of the cube $I^n$ (where $I=[0,1]$ is the interval) which we denote by $I^n_{t}$. Then the associahedron $\KKK^n$ is the quotient
$$\KKK^n:= \bigsqcup_{t}I^n_{t}/ \sim$$
where the equivalence relation is as follows. We think of an element $\tau=(t; \lambda_{1},\ldots ,\lambda_{n})\in I^n_{t}$ as a tree of type $t$ where the $\lambda_{i}$'s are the lengths of the internal edges. If some of the $\lambda_{i}$'s are 0, then the geometric tree determined by $\tau$ is not binary anymore (since some of its internal edges have been shrinked to a point). We denote the new tree by $\bar{\tau}$. For instance, if none of the $\lambda_{i}$'s is zero, then $\bar{\tau}=t$ ; if all the $\lambda_i$'s are zero, then the tree $\bar{\tau}$ is the corolla (only one vertex). The equivalence relation $\tau \sim \tau'$ is defined by the following two conditions:

-  $\bar{\tau}=\bar{\tau'}$,

-  the lengths of the nonzero-length edges of $\tau$ are the same as those of $\tau'$. 

 Hence $\KKK^n$ is obtained as a cubical realization denoted $\Kcub n$.
 
 Examples:
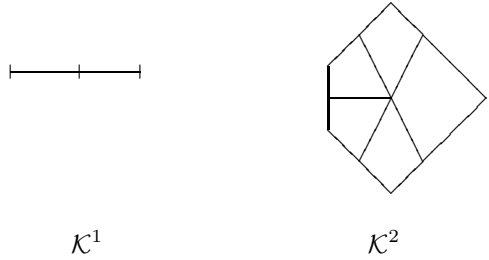
$$ \quad \KunSquare   \quad \KdeuxSquare$$
$$\KKK^1\qquad \qquad \qquad \qquad \qquad   \KKK^2 $$

\subsection{Markl-Shnider version of Saneblidze-Umble diagonal \cite{MS, SU1}} In \cite{BV} Boardman and Vogt showed that the bar-cobar construction on the operad $As$ is a dgns operad $\Omega B As$ whose chain complex in arity $n$ can be identified with the chain complex of the cubical decomposition of the associahedron:
$$(\Omega B As)_{n}= C_{*}(\Kcub {n-2})\ .$$
In \cite{MS} (where $\Kcub {n-2}$ is denoted $W_{n}$ and $\KKK^{n-2}$ is denoted $K_{n}$) Markl and Shnider use this result to construct a coassociative diagonal on the operad $\Omega B As$. There is a  quasi-isomorphism $q:\Ai \to \Omega B As$ induced by the cubical decomposition of the associahedron (the image of the top cell is the the algebraic sum of the $c_{n-1}$ cubes). They construct an inverse quasi-isomorphism $p:\Omega B As\to \Ai$ by giving explicit algebraic formulas. 
At the chain level the map $p: C_{*}(\Kcub {n}) \to C_{*}(\KKK^n)$  has a topological interpretation using a deformation of the cubical associahedron as follows. The cube indexed by the left comb is called the \emph{main cube} of the decomposition. The deformation sends the main cube to the top cell of the associahedron and flatten all the other ones. 

Example:

$$\KdeuxSquareMain  \KdeuxCubical$$

The exact way the main cube is deformed is best explained by drawing the associahedron on the cube. This is recalled in the Appendix. In \cite{KS} Kadeishvili and Saneblidze give a general method for constructing a diagonal on some polytopes admitting a cubical decomposition along the same principle (inflating the main cube).

Markl and Shnider claim that the composite 
$$ \Ai \stackrel{q}{\to} \Omega B As \to \Omega B As \t \Omega B As \stackrel{p\t p}{\longrightarrow} \Ai \t \Ai$$
is the Saneblidze-Umble diagonal.
\subsection{Comparison of the operads $\AAi$ and $\Omega B As$} From the geometric nature of $\AAi$ and $\Omega B As$ we compare them as follows. A cube admits a simplicial decomposition. Hence we can simplicialize the cubical decomposition of $\Kcub {n-2}$ to obtain a new cellular complex $\KKK_{\rm cub,simp}^{n}$. There are explicit quasi-isomorphisms
$$C_{*}(\Kcub {n})\to C_{*}(\KKK_{\rm cub,simp}^{n})\leftarrow C_{*}(\Ksimp n)\ .$$
We leave it to the reader to figure out the explicit formulas from the 2-dimensional case:

$$\KcubesimpUn \KcubesimpDeuxbis  \KdeuxT$$

Observe that the intermediate spaces $\KKK_{\rm cub,simp}^{n}$ give rise to a new dgns operad along the same lines as before.

\section {Appendix 1: Drawing a Stasheff polytope on a cube}

This is an account of some effort to construct the Stasheff polytope that I did in 2002 while visiting Northwestern University. During this visit I had the opportunity to meet Samson Saneblidze and Ron Umble, who were drawing the same kind of figures for different reasons (explained above).
It makes the link between Markl and Shnider algebraic description of the map $p$, the pictures appearing in Saneblidze and Umble paper, and some algebraic properties of the planar binary trees.

There is a way of constructing an associahedron structure on a cube as follows. For $n=0$ and $n=1$ there is nothing to do since $\KKK^0$ and $\KKK^1$ are the cubes $I^0$ and $I^1$ respectively. For $n=2$, we simply add one point in the middle of an edge to obtain a pentagon:
$$
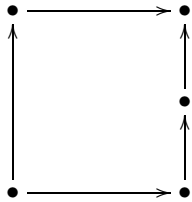
\xymatrix{
\bullet\ar[rr] & & \bullet \\
 & & \bullet\ar[u] \\
 \bullet\ar[uu]\ar[rr] & & \bullet\ar[u]\\
 }$$
 
 Inductively we draw $\KKK^n$ on $I^n$ out of the drawing of $\KKK^{n-1}$ on $I^{n-1}$ as follows. Any tree $t\in PBT_{n+1}$ gives rise to an ordered sequence of trees $(t_{1}, \ldots , t_{k})$ in $PBT_{n+2}$ as follows. We consider the edges which are on the right side of $t$, including the root. The tree $t_{1}$ is formed by adding a leaf which starts from the middle of the root and goes rightward (see \cite{JLLarithmetree} p.~297). 
 The tree $t_{2}$ is formed by adding a leaf which starts from the middle of the next edge and goes rightward. And so forth. Obviously $k$ is the number of vertices lying on the right side of $t$ plus one (so it is always greater than or equal to 2).
 
 Example:
 
  if $t= \arbreBA$, then $t_{1}= \arbreBAC, t_{2}= \arbreCAB, t_{3}= \arbreCBA$.
 
 In $I^n=I^{n-1}\times I$ we label the point $\{t\}\times\{0\}$ by $t_{1}$, the point $\{t\}\times\{1\}$ by $t_{k}$, and we introduce (in order) the points $t_{2}, \ldots , t_{k-1}$ on the edge $\{t\}\times I$. For $n=2$ we obtain (with the coding introduced in section \ref{planarbinarytree}):
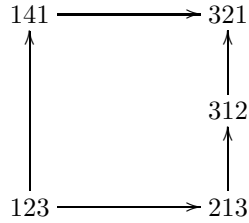
 $$\xymatrix{
141\ar[rr] & &321 \\
 & & 312\ar[u] \\
123\ar[uu]\ar[rr] & & 213\ar[u]\\
 }$$
 
 For $n=3$ we obtain the following picture:
 
 $$\KtroisCube$$
 
 (It is a good exercise to draw the tree at each vertex). Compare with \cite{SU1}, p.~3). The case $n=4$ can be found on my home-page. It is important to observe that the order induced on the vertices by the canonical orientation of the cube coincides precisely with the Tamari poset structure.
 
 Surprisingly, this way of viewing the associahedron is related to an algebraic structure on the set of planar binary trees $PBT = \bigcup_{n\ge 1} PBT_{n}$, related to dendriform algebras. Indeed there is a non-commutative monoid structure on the set of homogeneous nonempty subsets of $PBT$ constructed in \cite{JLLarithmetree}. It comes from the associative structure of the free dendriform algebra on one generator. This monoid structure is denoted by $+$, the neutral element is the tree $|\ $. If $t\in PBT_{p} $ and $s\in PBT_{q}$, then $s+t$ is a subset of $PBT_{p+q-1}$. It is proved in \cite{JLLarithmetree} that the trees which lie  on the edge $\{t\}\times I\subset I^n$ are precisely the trees of $t+\arbreA$. For instance:
$$ \arbreAB + \arbreA = \arbreABC \cup \arbreACA$$
and
$$ \arbreBA + \arbreA = \arbreBAC \cup \arbreCAB\cup \arbreCBA\ .$$

The deformation of the associahedron consisting in inflating the main simplex to the top cell can be performed into two steps by considering a cube inside the associahedron. This cube is determined by the previous construction. First, we inflate the main simplex to the full cube as described in \ref{deformedCube}, then we deform the cube into the associahedron as indicated above.

Finally we remark that the deformation described in  \ref{deformedAssociahedron} permits us to draw the associahedron on the simplex.

\section {Appendix 2: $\DD(\mu_{5})$} In this appendix we give the computation of 
$\DD(\mu_{5})$ and we show that we get the same result as Saneblidze and Umble. In order to compare with their result we adopt their way of indexing the planar trees, which is as follows. Let $t$ be a tree whose root vertex has $k+1$ inputs, that we label (from left to right) by $0, \ldots , k$. Then, by definition, $d_{ij}(t)$ is the tree obtained by replacing, locally, the root vertex by the following tree with one internal edge:

$$
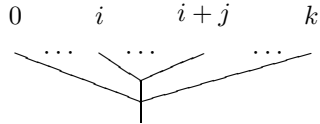
\xymatrix@R=2pt@C=2pt{
0 & & i & & i+j & & k \\
*{}\ar@{-}[ddrrr] & \cdots & *{}\ar@{-}[dr] & \cdots &  *{}\ar@{-}[dl] & \cdots & *{}\ar@{-}[ddlll] \\
& & & *{}\ar@{-}[dd] & & & \\
& & & *{}& & & \\
& & & *{}& & & 
}$$

The operator $d_{ij}$ is well-defined for $0\leq i\leq k, 1\leq j\leq k-i$ and $(i,j)\neq (0,k)$. So we get:
$$
\begin{array}{|cccccccc|}
\hline 
i j   &=&				& 01 & 02 & 11 & 12 & 21 \\
\hline
d_{ij}\Big(\arbreCCC\Big)&=&  & \arbreACC  & \arbreBBC   & \arbreCAC   & \arbreCBB   & \arbreCCA\\
\hline
\end{array}$$
and $d_{01}d_{01}\Big(\arbreCCC\Big)=\arbreABC$, etc.

Let us index the sixteen $3$-simplices forming $\Ksimp 3$ by the tree indexing the face in $fat\Ksimp 2$ and either $a,b,c$ if this face is a pentagon (cf.~\ref{simplicialSP}) or the shuffle $ \aa=(1,2), \bb=(2,1)$ if this face is a square (cf.~\ref{AAinfini}). In the following tableau we indicate the image of the $3$-simplices under the map $p'\t p' \circ \DD^{\AAi}$. In the left column we indicate the information which determines the $3$-simplex  $(d_{ij}(\mu_{5}), x)$. In the right column we give its image (up to signs) as a sum of four terms, since in the AW morphism there are four terms.
$$
\begin{array}{|c|c|}
\hline
03\quad a      &(01)(01)(01)\t \mu_{5}+  (02)(01) \otimes \big((21)+ (22)\big) \\
 & + (03)\t \big((11)(21)+(12)(21)+(11)(22) \big) +\mu_{5} \t (11)(21)(31)\\
03\quad b     &0+0+0+0 \\
03\quad c      &0+0+0+0 \\
02\quad \aa&0+0+(02) \t \big( (11)(31)+ (12)(31)\big)+0 \\
02\quad \bb&   0+ (01)(02) \t \big( (11) + (12) + (13) \big) +0+0 \\
01\quad a     & 0+ (01)(01) \t (31) + (01)\t (21)(31) \\
01\quad b     &  0+0+0+0 \\
01\quad c      &  0+0+0+0 \\
12\quad \aa&  0+0+(12)\t \big( (12)(21)+ (11)(22)\big) +0 \\
12\quad \bb&  0+0+0+0 \\
11\quad a      & 0+0+(11) \t (12)(31)+0 \\
11\quad b      &0+(02)(11)\t \big((13)+(12) \big)+0+0 \\
11\quad c      &  0+(11)(11)\t (13)+0+0 \\
21\quad a     &  0+(11)(01)\t (22)+(21)\t (11)(22)+0 \\
21\quad b     &  0+0+0+0 \\
21\quad c     & 0+0+0+0 \\
\hline
\end{array}$$

As a result $\DD(\mu_{5})$ is the algebraic sum of 22 elements, which are exactly the same as in \cite{SU1} Example 1. Topologically, it means that $\KKK^3$ can be realized as the union of 2 copies of $\KKK^3$ (having only one vertex in common), 6 copies of $\KKK^1\times \KKK^2$,  6 copies of $\KKK^2\times \KKK^1$,  4 copies of $(\KKK^1\times \KKK^1)\times \KKK^1 $ and   4 copies of $\KKK^1\times (\KKK^1\times \KKK^1)$.

From this computation it is reasonable to conjecture that  the diagonal constructed from the simplicial decomposition of the associahedron is the same as the Saneblidze-Umble diagonal.


\begin{thebibliography}{aa}

\bibitem{BV} \emph{J.M.~Boardman, R.M.~Vogt}, Homotopy invariant algebraic structures on topological spaces. Lecture Notes in Mathematics, Vol. 347. Springer-Verlag, Berlin-New York, 1973. x+257 pp.

\bibitem{GZ} \emph{M.~Gaberdiel, B.~Zwiebach}, 
Tensor constructions of open string theories. I. Foundations.
Nuclear Phys. B 505 (1997), no. 3, 569--624. 

 \bibitem{KS} \emph{T.~Kadeishvili, S.~Saneblidze}, The
twisted Cartesian model for the double path fibration, ArXiv math.AT/0210224

\bibitem{Keller}\emph{B.~Keller},
Introduction to $A$-infinity algebras and modules.
Homology Homotopy Appl. 3 (2001), no. 1, 1--35.

\bibitem{JLLarithmetree} \emph{J.-L.~Loday},  
Arithmetree.
J. Algebra 258 (2002), no. 1, 275--309.

\bibitem{JLLStasheff} \emph{J.-L.~Loday},  Realization of the Stasheff polytope. Arch. Math. (Basel) 83 (2004), no. 3, 267--278.

\bibitem{JLLparking} \emph{J.-L.~Loday}, Parking functions and triangulation of the associahedron, Proceedings of the Street's fest, Contemporary Math. AMS 431, (2007), 327--340.

\bibitem{MacLane} \emph{S.~MacLane}, {\bf Homology}. Die Grundlehren der mathematischen Wissenschaften, Bd. 114 Academic Press, Inc., Publishers, New York; Springer-Verlag, Berlin-G\"ottingen-Heidelberg 1963 x+422 pp.

\bibitem{MSS} \emph{M.~Markl, S.~Shnider, J.~Stasheff}, {\bf Operads in algebra, topology and physics}. Mathematical Surveys and Monographs, 96. American Mathematical Society, Providence, RI, 2002. x+349 pp.

\bibitem{MS} \emph{M.~Markl, S.~Shnider}, Associahedra, cellular $W$-construction and products of $A\sb \infty$-algebras. Trans. Amer. Math. Soc. 358 (2006), no. 6, 2353--2372 (electronic).

\bibitem{Proute} \emph{A.~Prout\'e}, $A_{\infty}$-structures, mod\`ele minimal de Baues-Lemaire et homologie des fibrations, Th\`ese d'Etat, 1984, Universit\'e Paris VII. 

\bibitem{SU1} \emph{S.~Saneblidze, R.~Umble},   A Diagonal on the Associahedra, preprint, ArXiv  math.AT/0011065 

\bibitem{SU2} \emph{S.~Saneblidze, R.~Umble},  Diagonals on the permutahedra, multiplihedra and associahedra. Homology Homotopy Appl. 6 (2004), no. 1, 363--411.

\bibitem{Stasheff} \emph{J.D. ~Stasheff}, Homotopy associativity of $H$-spaces. I, II. Trans. Amer. Math. Soc. 108 (1963), 275-292; ibid. 293--312. 



\end{thebibliography}
\end{document}